\documentstyle[12pt,amscd,amssymb]{amsart}

\newcommand{\hH}{\hat{H}}
\newcommand{\Dws}{D^{(w,\Sigma)}}
\newcommand{\M}{M^{\text{\scriptsize odd}}_{\Sigma}}
\newcommand{\Y}{\Sigma \times {{\Bbb S}}^1}
\newcommand{\X}{X^o}
\newcommand{\D}{D^o}
\newcommand{\B}{B^o}

\newcommand{\inc}{\hookrightarrow}
\newcommand{\ar}{\rightarrow}
\newcommand{\bd}{\partial}
\newcommand{\x}{\times}
\newcommand{\iso}{\cong}
\newcommand{\isom}{\stackrel{\sim}{\ar}}

\newcommand{\im}{\hbox{im}}
\newcommand{\pt}{\text{pt}}
\newcommand{\CP}{{\Bbb C \Bbb P}}

\newcommand{\Diff}{\text{Diff}}

\newcommand{\Sym}[1]{\text{Sym}^{#1}}

\newcommand{\red}{\hbox{\scriptsize red}}
\newcommand{\PD}{\text{P.D.}}

\newcommand{\cF}{{\cal F}}

\newcommand{\cH}{{\cal H}}

\renewcommand{\AA}{{\Bbb A}}

\newcommand{\DD}{{\Bbb D}}

\newcommand{\RR}{{\Bbb R}}
\newcommand{\SS}{{\Bbb S}}
\newcommand{\ZZ}{{\Bbb Z}}

\renewcommand{\a}{\alpha}
\renewcommand{\b}{\beta}

\newcommand{\g}{\gamma}

\newcommand{\f}{\epsilon}

\newcommand{\k}{\kappa}

\newcommand{\p}{\phi}
\newcommand{\q}{\psi}

\newcommand{\z}{\zeta}

\renewcommand{\S}{\Sigma}
\newcommand{\De}{\Delta}

\theoremstyle{plain}
\begingroup
\newtheorem{thm}{Theorem}
\newtheorem{cor}[thm]{Corollary}
\newtheorem{lem}[thm]{Lemma}
\newtheorem{prop}[thm]{Proposition}
\endgroup

\theoremstyle{definition}
\newtheorem{defn}[thm]{Definition}

\theoremstyle{remark}
\newtheorem{rem}[thm]{Remark}

\title[Donaldson invariants for connected sums]{Donaldson 
invariants for connected sums along surfaces of genus $2$}

\author{Vicente Mu\~noz}
\address{Departamento de \'Albegra, Geometr\'{\i}a y Topolog\'{\i}a  \\
Facultad de Ciencias \\
Universidad de M\'alaga, AP.\ 59 \\ 29080 M\'alaga \\ Spain}
\email{vmunoz@@agt.cie.uma.es}
\date{December, 1996}

\begin{document}

\maketitle

\begin{abstract}
  We prove a gluing formula for the Donaldson invariants of the connected
  sum of two four-manifolds along a surface of genus $2$. We also prove
  a finite type condition for manifolds containing a surface of genus
  $2$,
  self-intersection zero and representing an odd homology class.
\end{abstract}

\section{Introduction}
  \label{sec:intro}
  This paper tries to answer the question of the behaviour of the
  Donaldson
  invariants under connected sums along surfaces of genus $2$. This has
  been
  treated by the author in~\cite{Thesis} making use of a suitable version 
  of the Atiyah-Floer Conjecture. The purpose of this paper is to remove
  the use of any conjecture as well as to make the argument direct and 
  very simple. This was prompted by Tom Mrowka in the conference on
  four-manifolds in Oberwolfach (Germany) on May 96.
  We also remark that similar cases have been treated by Morgan and
  Szab\'o~\cite{Szabo}~\cite{Szabo2}, but our results are more general.

  Let $X$ be a smooth, compact, oriented four-manifold with $b^+ >1$
  and $b^+ -b_1$ odd. For 
  any $w \in H^2(X;\ZZ)$, $D^w_X$ will denote the corresponding Donaldson
  invariant~\cite{KM}, which is defined as a linear functional on 
  ${\AA(X)}= \Sym{*}(H_0(X)\oplus H_2(X))$
  ($H_*(X)$ will always denote homology with rational coefficients, and
  similarly for $H^*(X)$). Let $x \in H_0(X)$ be the class of a point.
  Then Kronheimer and Mrowka~\cite{KM} define $X$ to be of simple type 
  (with respect to $w$) when $D_X^{w}((x^2 -4)z)=0$ for
  all $z \in {\AA(X)}$, and in that case define 
$$ 
  \DD_X^{w}(z)=D_X^{w}((1+{x \over 2})z),
$$
  for all $z \in \Sym{*}H_2(X)$. The series $\DD_X^{w}(e^{t\a})$, $\a \in
  H_2(X)$, is even or odd depending on whether
  $d_0=d_0(X,w)= -w^2- {3\over 2}(1-b_1+b^+)$ is even or odd.
  When $b_1=0$ and $b^+>1$, $X$ is of simple
  type with respect to some $w$ if and only if it is so with respect to
  any $w$. In such case, $X$ is just called of simple type.

\begin{prop}
\label{prop:KM}
  Let $X$ be a manifold of simple type with $b_1=0$ and 
  $b^+ >1$ and odd. Then we have 
  $$
     \DD_X^{w}(e^{\a})= e^{Q(\a) /2} \sum (-1)^{{K_i \cdot {w} +{w}^2}
     \over 2} a_i \, e^{K_i \cdot \a}
  $$
  for finitely many $K_i \in H^2(X;\ZZ)$ (called basic
  classes) and rational numbers $a_i$ (the collection is empty
  when the invariants all vanish). These classes are lifts
  to integral cohomology of ${w}_2(X)$. Moreover, for any embedded
  surface
  $ S \inc X$ of genus $g$ and with $S^2 \geq 0$,
  one has $2g-2 \geq S^2 +|K_i \cdot S|$. 
\end{prop}

  Analogously, we define $X$ to be of finite type with 
  respect to $w$ whenever $D_X^{w}((x^2 -4)^n z)=0$ for
  all $z \in {\AA(X)}$, and some $n>0$. The order is the minimum
  of such $n$.

  When $X$ has $b^+=1$, the invariants depend on the metric through a 
  structure of walls and chambers~\cite{Kotschick}
  and therefore we have to specify the metric.

\begin{defn} 
\label{def:allowable}
  $(w,\S)$ is an {\bf allowable} pair if 
  $w, \S \in H^2(X; \ZZ)$, $w \cdot \S \equiv 1\pmod 2$ and $\S^2 =0$. 
  Then we define
  $$
    D^{(w,\S)}_X=D^w_X +D^{w+\S}_X.
  $$
  When $b^+=1$ we consider the invariants referring to
  the chambers defined by $\S$, i.e. for metrics
  whose period points are in the (unique) chamber containing $\S$ in its
  closure (which is so since $w \cdot \S \equiv 1\pmod 2$). In fact, we
  would need a result saying that the invariants only depend on the 
  metric through the period point. This is true for simply-connected 
  manifolds and for $\S\x \CP^1$, with $\S$ a Riemann surface, which are
  all the cases we need for our arguments.
\end{defn}

  $D^{(w,\S)}_X$ depends only on $\S$
  and $w\pmod{\S}$, since $D^{w+2\S}_X=D^w_X$.
  As $(w+\S)^2 \equiv w^2 +2 \pmod 4$, we can recover $D^w_X$ and
  $D^{w+\S}_X$ from
  $D^{(w,\S)}_X$. The series $D^{(w,\S)}_X(e^{t\a})$, $\a \in
  H_2(X)$, is even or odd according to whether 
  $d_0$ is even or odd.

\begin{prop}
\label{prop:DwS}
  Suppose $X$ is a manifold of simple type with $b_1=0$ and $b^+>1$
  and odd. Write the Donaldson series as $\DD_X^w (e^{\a})=
  e^{Q(\a)/2} \sum (-1)^{{K_j \cdot {w} +{w}^2} \over 2}a_j\,
  e^{K_j \cdot \a}$. 
  Then setting $d_0=d_0(X,w)=-w^2-{3\over 2}(1+b^+)$ we have 
  $$
    D^{(w,\S)}_X(e^{\a})=e^{Q({\a})/2} \hspace{-8mm} \sum_{K_j \cdot
    \S \equiv 2\pmod 4} \hspace{-8mm} (-1)^{{K_j \cdot {w} +{w}^2}
    \over 2} a_j e^{K_j \cdot \a} + e^{-Q(\a)/2} \hspace{-8mm} \sum_{K_j
    \cdot \S \equiv 0\pmod 4} \hspace{-8mm} i^{-d_0} (-1)^{{K_j \cdot
    {w} +{w}^2} \over 2} a_j e^{iK_j \cdot \a}
  $$
  So giving $\DD_X^{\,w}$ is equivalent to giving $D^{(w,\S)}_X$. 
\end{prop}

\begin{pf}
  Note that $K_j \cdot \S \equiv 0 \pmod 2$, for all basic classes $K_j$.
  Since $((w+\S)^2+K_j \cdot (w+\S)) = (w^2 + K_j \cdot w) + 2 (w \cdot
  \S + K_j \cdot \S /2)$ we have
  $$
    \DD_X^{w+\S} (e^{\a})= e^{Q(\a)/2} \hspace{-8mm} \sum_{K_j \cdot \S
    \equiv 2\pmod 4} \hspace{-8mm} (-1)^{{K_j \cdot {w} +{w}^2} \over
    2}a_j\, e^{K_j \cdot \a} - e^{Q(\a)/2} \hspace{-8mm} 
    \sum_{K_j \cdot \S
    \equiv 0\pmod 4} \hspace{-8mm} (-1)^{{K_j \cdot {w} +{w}^2} \over
    2}a_j \, e^{K_j \cdot \a}
  $$
  Now since the only powers in $D_X^w(e^{t\a})$ are those $t^d$ with
  $d \equiv d_0\pmod 4$ one has
  $$ 
     D_X^w(e^{t\a})= {1 \over 2}(\DD_X^w(e^{t\a}) +
     i^{-d_0}\DD_X^w(e^{it\a}))
  $$ 
  and analogously 
  $$ 
     D_X^{w+\S}(e^{t\a})= {1 \over 2}(\DD_X^{w+\S}(e^{t\a}) -
     i^{-d_0}\DD_X^{w+\S} (e^{it\a}))
  $$
  since $d_0(X,w+\S)=d_0(X,w)+2$. So we finally get
  $$
    D^{(w,\S)}_X(e^{\a})=e^{Q({\a})/2} \hspace {-8mm} \sum_{K_j \cdot
    \S \equiv 2\pmod 4} \hspace{-8mm} (-1)^{{K_j \cdot {w} +{w}^2}
    \over 2} a_j e^{K_j \cdot \a} + i^{-d_0} e^{-Q(\a)/2} \hspace{-8mm}
    \sum_{K_j
    \cdot \S \equiv 0\pmod 4} \hspace{-8mm} (-1)^{{K_j \cdot
    {w} +{w}^2} \over 2} a_j  e^{iK_j \cdot \a}
  $$
\end{pf}

\begin{defn}
\label{def:permissible}
  We say that $(X,\S)$ is {\bf permissible} if $X$ is a smooth compact
  oriented four-manifold and $\S \inc X$ is an embedded Riemann surface
  of genus $2$ and self-intersection zero such that $[\S] \in H_2(X;\ZZ)$
  is odd (its reduction modulo $2$ is non-zero, or equivalently, it is an
  odd multiple of a primitive homology class). So we can consider
  $w \in H^2(X;\ZZ)$ with $w \cdot \S \equiv 1 \pmod 2$. Then $(w,\S)$ 
  is an allowable pair. This implies that $b^+>0$. Let $N_{\S} \iso
  A =\S \x D^2$ be an open tubular neighbourhood of $\S$ and 
  set $\X=X - N_{\S}$. Then 
  $\bd \X=Y \iso \Y$ (but the isomorphism is not canonical). We consider
  one such isomorphism fixed and (when necessary) we furnish $\X$ with a
  cylindrical end, i.e. we consider $\X \cup (Y\x [0,\infty))$ (and keep
  on calling it $\X$).
\end{defn}

  We call {\bf identification} for $Y=\Y$ any (orientation preserving) 
  bundle automorphism $\p: Y \isom Y$. Up to isotopy, $\p$ depends only 
  on the isotopy class of the induced diffeomorphism on $\S$ and on an
  element of $H^1(\S;\ZZ)$.

\begin{defn}
\label{def:sum}
  Let $(X_1,\S_1)$ and $(X_2,\S_2)$ be permissible. We pick orientations
  so that $\bd \X_1 = -\bd \X_2=Y$ (minus means reversed orientation).
  Then $X=X(\p) = \X_1 \cup_{\p} \X_2 = X_1 \#_{\S} X_2$ is a compact,
  naturally oriented, smooth four-manifold, called the {\bf connected 
  sum along $\S$} of $(X_1,\S_1)$ and $(X_2,\S_2)$ (with identification
  $\p$). The induced homology classes $[\S_1]$ and $[\S_2]$ coincide and
  are induced by an embedded $\S \inc X$. Then $(X,\S)$ is permissible.
\end{defn}

  Choose $w_i \in H^2(X_i;\ZZ)$, $i=1,2$, and $w \in H^2(X;\ZZ)$ such 
  that $w_i \cdot \S_i \equiv 1 \pmod 2$, $w \cdot \S \equiv 1 \pmod 2$, 
  in a compatible way (i.e. the restricition of $w$ to $\X_i \subset X$ 
  coincides with the restriction of $w_i$ to $\X_i \subset X_i$).
  We shall call $w$ all of them, not making explicit to which manifold
  they refer. Also let $\cH=\{ D \in H_2(X)/ D|_Y=k [\SS^1] \in H_1(Y),
  \hbox{ for some $k$}\}$. Then for every $D \in \cH$, we can
  choose $D_i \in H_2(X_i)$ agreeing
  with $D$ (i.e. $D_i|_{\X_i}=D|_{\X_i}$, $i=1,2$) 
  and with $D^2 =D_1^2 +D_2^2$. Furthermore, we can arrange
  $D \mapsto (D_1, D_2)$ to be linear. Once chosen one of these maps,
  any other is of the form $D \mapsto (D_1-r\S, D_2+ r\S)$, for a
  rational number $r$.
  
  A simple but important remark is that if $b_1(X_1)=b_1(X_2)=0$ then
  $b_1(X)=0$ and $b^+(X)>1$. Now we are ready to state our main results.

\begin{thm}
\label{thm:main1}
  Suppose $(X_1,\S_1)$ and $(X_2,\S_2)$ are permissible and $X_1$, $X_2$ 
  have both $b_1=0$ and $b^+>1$ and are of simple type. Let\/
  $\DD^w_{X_1}(e^{\a})=e^{Q(\a)/2}\sum a_{i,w} e^{K_i\cdot \a}$ and\/
  $\DD^w_{X_2}(e^{\a})=e^{Q(\a)/2}\sum b_{j,w} e^{L_j\cdot \a}$.
  Let $X=X_1 \#_{\S} X_2$ (for some identification). 
  Then $X$ is of simple type and
  for every $D \in \cH$, choose $D_i \in H_2(X_i)$ agreeing
  with $D$ satisfying 
  $D^2 =D_1^2 +D_2^2$, in such a way that $D \mapsto (D_1, D_2)$
  is linear. Then
  $$\DD^w_X(e^{tD}) =$$
  $$ = e^{Q(tD)/2}(\hspace{-5mm}\sum_{K_i \cdot \S=L_j \cdot
    \S = 2}\hspace{-5mm} - 32a_{i,w}b_{j,w} \, e^{(K_i \cdot D_1 +  L_j
    \cdot D_2 +2\S\cdot D)t} +\hspace{-5mm}\sum_{K_i \cdot
    \S=L_j \cdot \S= -2}\hspace{-5mm} 32a_{i,w}b_{j,w} \, e^{(K_i \cdot
    D_1 +  L_j \cdot D_2-2\S\cdot D)t}),
  $$
  (for appropriate homology orientations).
\end{thm}

\begin{rem}
\label{rem:signs}
  The reason for the different signs is easy to work out. First, $w^2$ 
  for $X$ is always congruent $\pmod 2$ with the sum of both of $w^2$ 
  for $X_i$. Also $-{3 \over 2}(1-b_1(X)+b^+(X))= -{3 \over
  2}(1-b_1(X_1)+b^+(X_1))-{3 \over 2}(1-b_1(X_2)+b^+(X_2)) -3(g-1)$.
  Therefore, as $g=2$, $d_0(X,w) \equiv d_0(X_1,w)+d_0(X_2,w)+1 \pmod 2$. 
  Now the sign comes from the fact that
  the coefficient for the basic class $-\k$ is $(-1)^{d_0}c_{\k}$, being
  $c_{\k}$ the coefficient for the basic class $\k$.
\end{rem}

\begin{cor}
\label{cor:Vic}
  Suppose we are in the conditions of the former theorem. 
  Write\/ $\DD_X(e^{\a})=e^{Q(\a)/2} \sum c_{\k}
  e^{\k\cdot \a}$ for the Donaldson series for $X$ and
  $\DD_{X_1}(e^{\a})=e^{Q(\a)/2}\sum a_i e^{K_i\cdot \a}$ and\/
  $\DD_{X_2}(e^{\a})=e^{Q(\a)/2}\sum b_j e^{L_j\cdot \a}$ for the Donaldson
  series for $X_1$ and $X_2$. 
  Then given any pair $(K,L) \in H^2(\X_1;\ZZ)
  \oplus H^2(\X_2;\ZZ)$, we have
  $$
   \sum_{\{\k/ \k|_{\X_1}=K, \, \k|_{\X_2}=L\}}\hspace{-4mm} c_{\k} = \pm
    32\,(\sum_{K_i|_{\X_1}=K} a_{i}) \cdot  (\sum_{L_j|_{\X_2}=L}
    b_{j})  
  $$
  whenever $K|_Y = L|_Y = \pm 2 \PD [\SS^1]$. Otherwise, the left hand
  side is zero.
\end{cor}

\begin{pf}
  This is an immediate consequence of theorem~\ref{thm:main1}, noting that
  $$(-1)^{\k \cdot w +w^2 \over 2}= - (-1)^{K_i \cdot w +w^2 \over 2}
  (-1)^{L_j \cdot w +w^2 \over 2}$$
  whenever $\k|_{\X_1}=K_i|_{\X_1}=K$, $\k|_{\X_2}=L_j|_{\X_2}=L$.
\end{pf}

\begin{thm}
\label{thm:main2}
  Let $(X,\S)$ be permissible. Then $X$ is of finite type (with respect
  to
  any $w \in H^2(X;\ZZ)$ with $w \cdot \S \equiv 1 \pmod 2$, and for the
  invariants given by $\S$ in case $b^+=1$).
\end{thm}

Now we introduce a very important example. Let $B$ be the $K3$-surface
blown-up in two points. Let $S \subset K3$ be a tight surface of genus $2$ 
(which existence is guaranteed by~\cite{KM})
and let $E_1$, $E_2$ be the two exceptional divisors in $B$. Then 
$\S=S-E_1-E_2$ is the proper transform of $S$. So $(B,\S)$ is permissible.
For any $(X, \S)$ permissible, write 
$\tilde X=X \#_{\S} B$ (fixing some identification). It has $b^+(X)>1$
and $b_1(\tilde X)=0$ whenever $b_1(X)=0$. Now for any embedded surface
$\D \subset \X$ with $\bd \D =\bd \X \cap  \D$, we 
can choose cappings $D=\D+ \D_B$ in $\tilde X$
(in general it is enough to suppose that $\D$ is a cycle and that, giving 
$\X$ a cylindrical end, $\D \cap (Y \x [0,\infty))= \g \x [0,\infty)$, with
$\g \subset Y$ an embedded curve).
Fix an embedded surface representing $E_1+E_2+\S$ and intersecting $\S$ 
transversely in two points, and let
$K^o$ be its restriction 
to $B^o$. Then we always impose $\D_B \cdot K^o =0$ (this pairing makes sense
as long as $\bd K^o$ and $\bd \D_B$ are disjoint).

\begin{thm}
\label{thm:main3}
  Let $(X_i,\S_i)$ be permissible with $b_1(X_i)=0$, $i=1,2$ (not
  necessarily of simple type).
  Consider $\tilde X_i =X_i
  \#_{\S} B$. Then $\tilde X_i$ are of simple type. Put\/ $\DD^w_{\tilde
  X_1}(e^{\a})=e^{Q(\a)/2}\sum \tilde a_{i,w} e^{\tilde K_i \cdot  \a}$ 
  and\/ $\DD^w_{\tilde
  X_2}(e^{\a})=e^{Q(\a)/2}\sum \tilde b_{j,w} e^{\tilde L_j \cdot \a}$.
  Let $X=X_1 \#_{\S} X_2$ (for some  identification). 
  Then $X$ is of simple type. For every $D \in H_2(X)$, consider any
  cappings $D_i \in H_2(\tilde X_i)$ with the condition above
  in such a way that $D \mapsto (D_1,D_2)$ is linear. Then
$$
  \DD^w_X(e^{tD})= e^{Q(tD)/2}   (\hspace{-5mm} \sum_{\tilde K_i
  \cdot \S=\tilde L_j 
  \cdot \S= 2} \hspace{-5mm} -{1\over 2}\tilde a_{i,w}\tilde b_{j,w}
  e^{t(\tilde K_i \cdot D_1 + \tilde L_j \cdot D_2)} 
  + \hspace{-5mm} \sum_{\tilde K_i \cdot \S=\tilde L_j
  \cdot \S= -2} \hspace{-5mm} {1\over 2}\tilde a_{i,w}\tilde b_{j,w}
  e^{t(\tilde K_i \cdot D_1 + \tilde L_j \cdot D_2)} ).
$$
\end{thm}

\begin{cor}
\label{cor:Vic2}
  Under the conditions of theorem~\ref{thm:main3}, $X$ has no basic
  classes $\k$ with $\k \cdot \S=0$.
\end{cor}

\noindent {\em Acknowledgements:\/} First I want to thank my D. Phil.\
supervisor Simon Donaldson, for his encouragement and help during
the last three years. Also I am very grateful to the Mathematics 
Department in Universidad de M\'alaga for their hospitatility 
and for letting me use their facilities, and 
to Banco de Espa\~na for financial support. The author wishes to thank
T. Mrowka for suggesting he remove the use of the Atiyah-Floer
conjecture from his earlier work.

\section{Applications}
\label{sec:applic}

Now we pass on to give some nice and simple applications of
theorem~\ref{thm:main3}. Probably, many results like the following can be
obtained in the same fashion. We only want to give some examples
to show its usefulness.

\begin{cor}
\label{cor:main1}
  Let $(X_1,\S_1)$ and $(X_2,\S_2)$ be permissible with $b_1(X_i)=0$.
  Let $\p$ and $\q$ be two different identifications for $Y=\Y$ and
  consider the two different connected sums along $\S$, $X(\p)$ and
  $X(\q)$. 
  Suppose that $\p_* =\q_* : H_1(Y) \ar H_1(Y)$.
  Then there is an (non-canonical) 
  isomorphism of vector spaces $H^2(X(\p)) \isom H^2(X(\q))$
  sending the basic classes of $X(\p)$ to those of $X(\q)$ such that
  the rational numbers attached to them coincide.
\end{cor}

\begin{pf}
  First we observe that we have a natural identification of the images
  $I_{\p}$ of $H_2(\X_1) \oplus H_2(\X_2) \ar H_2(X(\p))$ and $I_{\q}$
  of 
  $H_2(\X_1) \oplus H_2(\X_2) \ar H_2(X(\q))$ since the kernels
  coincide. Now consider a splitting $H_2(X_2(\p)) \iso \text{Im}(I_{\p}) 
  \oplus V$ with
  $V \isom H_1(Y)$. 
  Choose an integral basis $\{\a\}$ for $H_1(Y;\ZZ)$.
  For every $\a$ we have an element
  $D_{\a} \in H_2(X(\p))$ which can be split as $D_{\a}=\D_1+\D_2$, for
  $\D_i \subset \X_i$ with $\bd \D_1 =\g$, $-\bd \D_2 =
  \p(\g)$ and $\a=[\g]$. Now we leave $\D_1$ (and $D_1 \in
  H_2(\tilde X_1)$)
  fixed and modify $\D_2$ to glue it to $\D_1$ in $H_2(X(\q))$.
  Write $D_2=\D_2+\D_3 \in H_2(\tilde X_2)$. The loops $\p(\g)$
  and $\q(\g)$ are 
  homologous and hence there is homology $C=\SS^1 \x [0,1] \inc \S
  \subset \Y$ between them. Consider 
  $$
  (D')^o_3=\left[ \D_3 \cup_{\p(\g)} C \cup_{\q(\g)}
  (\q(\g) \x [0,\infty))\right]+n\S \subset \B
  $$
  $$
  (D')^o_2=\left[ \D_2 \cup_{\p(\g)}(- C) \cup_{\q(\g)}
  (-\q(\g) \x [0,\infty))\right]-n\S \subset \X_2
  $$
  where $n$ is chosen so that $(D')^o_3 \cdot K^o =0$. So $D_2'
  =(D')^o_2+(D')^o_3= D_2$.
  Consider $D'_{\a}=\D_1+ (D')^o_2 \in H_2(X(\q))$. The map $D_{\a}
  \mapsto D'_{\a}$ gives the required isomorphism $H^2(X(\p)) \isom
  H^2(X(\q))$.
\end{pf}

  This corollary says that although in principle $X(\p)$ and $X(\q)$
  might not be diffeomorphic (and probably in many cases this happens),
  they can not be distinguished by the number and coefficients of their
  basic classes. Still the polynomial invariants can differentiate
  both manifolds (maybe the intersection matrix of the basic classes
  could help). It would be desirable to find examples when this happens.
  The identifications  to try out could be Dehn twists along separating 
  curves in $\S$.

\begin{cor}
\label{cor:main2}
  Let $(X_1,\S_1)$ and $(X_2,\S_2)$ be permissible with $b_1(X_i)=0$.
  Let $\p$ and $\q$ be two different identifications for $Y=\Y$ and
  consider the two different connected sums along $\S$, $X(\p)$ and 
  $X(\q)$. 
  Suppose that $X(\p)$ has only two basic classes $\pm \k$. Then the
  same is true for $X(\q)$ and the coefficients coincide (up to sign).
  Also if the invariants of $X(\p)$ vanish (no basic classes), so do
  the invariants of $X(\q)$.
\end{cor}

\begin{pf}
We do the case of two basic classes. The other one is analogous.
Suppose $\p=\text{Id}$, put
$X=X(\p)$ and let $\pm \k$ be the two basic classes, with
$\k \cdot \S=2$. Let $c_{\k,w}$ be its coefficient.
We now want to prove that this implies that there is only one basic
class $\tilde K_i$ with $\tilde K_i \cdot \S=2$ and only one basic
class $\tilde L_j$ with 
$\tilde L_j \cdot \S=2$. The result follows from
this applying theorem~\ref{thm:main3}.

Suppose that we can find $S_i \in H_2(\tilde X_i)$ with $\a= S_1\cap [Y] =-
S_2 \cap [Y] \in H_1(Y;\ZZ)$ such that all the values $\tilde K_i \cdot
S_1$ are different among them, and all the values $\tilde L_j \cdot
S_2$ are also different among them (where $\tilde K_i$ and $\tilde 
L_j$ run through all the basic
classes in $\tilde X_1$ and $\tilde X_2$ evaluating $2$ on $\S$). 
Then reorder the subindices in such a way that 
$$
  \tilde K_1 \cdot S_1 <\tilde K_2 \cdot S_1 <\cdots <\tilde 
    K_{n_1}\cdot S_1
$$  
$$
  \tilde L_1 \cdot S_2 <\tilde L_2 \cdot S_2 <\cdots <\tilde 
    L_{n_2}\cdot S_2
$$  
We can easily arrange $\D_i \subset \X_i$ with $\bd \D_1 = -\bd
\D_2=\g$ with $[\g]=\a$ such that the corresponding $D_i$ is $S_i$. Set
$D=\D_1 +\D_2 \in H_2(X)$ and apply theorem~\ref{thm:main3}. We have
$$
  c_{\k,w} e^{t\k \cdot D}  
  =\sum\limits_{\tilde K_i \cdot \S=\tilde L_j
  \cdot \S= 2} -{1 \over 2} \tilde a_{i,w}\tilde b_{j,w}
  e^{t(\tilde K_i \cdot S_1 + \tilde L_j \cdot S_2)} 
$$
Considering the exponentials with the smallest and with the largest
exponents, we see that we must have $\tilde K_1 \cdot S_1 + \tilde L_1
\cdot S_2 =\tilde K_{n_1} \cdot S_1+\tilde L_{n_2}\cdot S_2$, from
where the result.

To find the required collection of $S_i$, 
we consider all the differences $\a_{ij}=
\tilde K_i-\tilde K_j$, $\b_{ij}=\tilde L_i-\tilde L_j$, $i\neq j$. 
Consider $\a \in H_1(Y;\ZZ)$ such that $\a \cdot \a_{ij}\neq 0$ for any
$\a_{ij}$ which happens to be 
in the image of the homomorphism $H^1(Y) \iso H_2(Y) \inc
H_2(\tilde X_1) \iso H^2(\tilde X_1)$, and $\a \cdot \b_{ij} \neq 0$
when $\b_{ij}$ is in the same condition with $\tilde X_2$ replacing
$\tilde X_1$. Now we can choose $S_1 \in H_2(\tilde X_1)$ with $S_1
\cap [Y]=\a$ such that $\a_{ij} \cdot S_1 \neq 0$ (indeed the bad set
is a finite union of hyperplanes). Analogously we
choose $S_2$.
\end{pf}

\section{Gluing theory}
\label{sec:gluing}

  Let $X=\X_1 \cup_Y \X_2$, $D \in H_2(X)$. Substitute $D$ by a rational
  multiple if necessary so that $D|_Y \in H_1(Y;\ZZ)$ and it is
  primitive.
  Represent $D$ by a cycle so $D= \D_1+\D_2$, $\D_i \subset \X_i$,
  $\bd \D_1=-\bd\D_2=\g$, with $\g \subset Y$ an embedded curve in $Y$
  (when $\X_1$ has a cylindrical end, we suppose
  $\D_1 \cap (Y\x [0,\infty)) = \g \x [0,\infty)$, and analogously for
  $\X_2$).

\begin{prop}[\cite{BD}\cite{Thesis}]
\label{prop:Fukaya}
  Suppose $w|_Y$ odd. Then there are Fukaya-Floer homology groups 
  $HFF_*(Y,\g)$ 
  graded mod $4$ such that $(\X_i, \D_i)$ define relative invariants
  $\p^{w_1}(\X_1,e^{t\D_1}) \in HFF_*(Y,\g)$, 
  $\p^{w_2}(\X_2,e^{t\D_2})\in HFF_*(-Y,-\g)$. There is a natural pairing
  such that 
  \begin{equation}
     D_X^{(w,\S)}(e^{tD})=
     <\p^{w_1}(\X_1,e^{t\D_1}),\p^{w_2}(\X_2,e^{t\D_2})>. \label{eqn:sym}
  \end{equation}
  When $b^+=1$, the invariants are calculated for a long neck, i.e. we
  refer to the invariants defined by $\S$.
\end{prop}

In our case $Y \iso (-Y)$. Also, as explained in~\cite{BD}, 
$HFF_*(Y,\g)$ is
the limit of a spectral sequence whose $E_3$-term is $HF_*(Y) \otimes
\hH_*(\CP^{\infty})$ (the hat means the natural completion of 
$H_*(\CP^{\infty})$), and $d_3$ is multiplication by $\mu(\g)$.

First, $HF_*(Y)= HF_*(\Y) \iso HF_*^{\mathrm{symp}}(\M) \iso
H_*(\M)$ as vector spaces (we are using rational coefficients), 
where $\M$ is the moduli space of odd degree rank two 
stable vector bundles on $\S$ (with the grading considered mod $4$)
(for the first isomorphism see~\cite{Dostoglou}, for the second 
see~\cite{Schwarz}). For $g=2$, these groups were computed by
Donaldson~\cite{Don1}, finding that $H_*(\M)$ has an even part of 
dimension $4$ and an odd part of dimension $4$ (in the even part the
intersection product is symmetric, in the odd part it is antisymmetric).

There is a conjecture asserting that multiplication by
$\mu(\g)$ is intertwined with 
quantum multiplication by $\mu(\g)$ (see~\cite{Don1}~\cite{Thesis}). 
In~\cite[chapter 5]{Thesis}, the author has studied the implications of
such a conjecture. Here we want to avoid it altogether.

Essentially we have two cases to deal with, $\g=\SS^1 \subset \Y=Y$ and
$\g \subset \S \subset \Y=Y$. 

\begin{itemize}

\item $\g=\SS^1 \subset \Y$. Now all the
differentials in the $E_3$ term of the spectral sequence are of the
form $H_{\mathrm{odd}}(\M) \ar H_{\mathrm{even}}(\M)$ and
$H_{\mathrm{even}}(\M) \ar H_{\mathrm{odd}}(\M)$. When the boundary
cycle is $\g=\SS^1$ and thus invariant under the action of the group
$\Diff (\S)$ on $Y=\Y$, the differentials commute with the
action of $\Diff (\S)$. As there are elements $\rho \in
\Diff (\S)$ acting as $-1$ on $H^1(\S)$, we have that $\rho$
acts as $-1$ on $H_{\mathrm{odd}}(\M)$ and as $1$ on
$H_{\mathrm{even}}(\M)$. Therefore the differentials are zero and the
spectral sequence degenerates in the third term. 
This implies that 
$HFF_{\mathrm{even}}(Y,\g) = V_4 [[t]]$, where $V_4=HF_{\mathrm{even}}(Y)$
has dimension $4$. 
The relative invariants will be
$\p^{w_1}(\X_1,e^{t\D_1}) \in V_4[[t]]$. We do not consider the odd
part since the pairing is antisymmetric in the odd part, but 
the expression~\eqref{eqn:sym} is symmetric.

\item $\g \subset \S \subset \Y$.
The $E_3$ term of the usual spectral sequence is $HF_*(Y) \otimes 
\hH_*(\CP^{\infty})$, with differencital $d_3$ given by 
$$
   \mu(\g):H_i(\M) \otimes H_j(\CP^{\infty}) \ar H_{i-3}(\M) \otimes
   H_{j+2}(\CP^{\infty}).
$$

\begin{lem}
The image of $d_3:HF_3 \ar HF_0$ is one-dimensional and the kernel 
of $d_3:HF_2 \ar
HF_3$ is one-dimensional.
\end{lem}

\begin{pf}
Let us see first that
$$
  \mu(\g): HF_*(Y) \ar HF_*(Y)
$$ 
is non-zero. We decompose $\S \x \CP^1 = A \cup_Y A$. From the definition of
$\mu(\g)$ (see~\cite{Don2}~\cite{Braam}~\cite{Thesis}), we have that for 
$X=\X_1 \cup_Y \X_2$, $z_i  \in \AA(\X_i)$,
$\p^w(\X_i,z_i) \in HF_*(Y)$ and $\b \in H_*(Y)$, it is 
$\p^w(\X_1,\b z_1)= \mu(\b)(\p^w(\X_1,z_1))$. 
Also  $\Dws_X(z_1z_2) $ 
$=<\p^w(\X_1,z_1), \p^w(\X_2,z_2)>$.
We have thus
$$ 
  \Dws_{\S \x \CP^1}(\g_1\g_2) =< \p^w(A,1),\mu(\g_1)\mu(\g_2)
  (\p^w(A,1))>.
$$
The invariant of the left hand side corresponds to the six-dimensional
moduli space. This is in fact $\M$. From~\cite{Don1}~\cite{Thaddeus} 
we know that this number is non-zero (actually  $\f_S(w)\, \g_1\cdot \g_2$,
with $\f_S(w)=(-1)^{K_S \cdot w+w^2 \over 2}$). 
Therefore  $\mu(\g) \neq 0$.

Under the intersection pairing, $HF_2 \iso (HF_0)^*$ and $HF_3 \iso 
(HF_3)^*$. also, $d_3: HF_3 \ar HF_0$ and $d_3: HF_2 \ar HF_3$ are dual
maps, so the dimensions of $\ker(d_3: HF_2 \ar HF_3)$ and 
$\im(d_3: HF_3 \ar HF_0)$ coincide. 
Since $d_3$ is non-zero, these dimensions are at least one.
They cannot be two because that would imply that $HFF_{\mathrm{even}}(Y,\g) =
HF_0 \oplus 0 \oplus 0 \oplus \cdots$ and hence
$$
\Dws_X(e^{tD})= <\p^{w}(\X_1,e^{t\D_1}),\p^{w}(\X_2,e^{t\D_2})> =0
$$
for any case in which $D=\D_1 +\D_2$, $\D_i  \subset \X_i$, with
$\bd \D_1= -\bd \D_2 =\g$. In particular, the invariants $\Dws_X$ would
vanish whenever $X=X_1 \#_{\S} X_2$ with $b_1(X_i)=0$, $i=1,2$. But this is 
impossible, as we will see examples in the proof of theorem~\ref{thm:main3}
when the invariants do not vanish (these examples are independent of the 
computation of $HFF_*(Y,\g)$ for $\g \subset \S \subset Y$).
\end{pf}

From this we write the even part of the $E_5$ term of the spectral sequence.
Set $HF_2^{\mathrm{red}}=\ker(d_3: HF_2 \ar HF_3)$, $HF_0^{\mathrm{red}}=
HF_0/\im(d_3: HF_3 \ar HF_0)$. The even part of the $E_5$ term is
$$
\begin{array}{lcccccccc}
  (E_5)_0 &=& HF_0 & \oplus \;\; 0\;\;  \oplus & HF_2^{\red} &
  \oplus \;\; 0\;\;   \oplus & HF_0^{\red} & \oplus & \cdots \\
  (E_5)_2 &=& HF_2^{\red} & \oplus \;\; 0\;\;  \oplus &
  HF_0^{\red} & \oplus 
  \;\; 0\;\;  \oplus & HF_2^{\red} & \oplus & \cdots
\end{array}
$$
The differential $d_5$ has to be zero (at least on the even part of $E_5$),
since otherwise we would have again that the invariants $\Dws_X$ vanish 
whenever $X=X_1 \#_{\S} X_2$ with $b_1(X_i)=0$, $i=1,2$.

Hence $HFF_{\mathrm{even}}(Y,\g)$ is equal to this $(E_5)_{\mathrm{even}}$.
Write $HF_0=\RR \oplus HF_0^{\red}$, where $\RR$ is the orthogonal complement
to $HF_2^{\red}$. Then $HFF_{\mathrm{even}}(Y,\g)=\RR \oplus V_2[[t]]$, with
$V_2= HF_0^{\red} \oplus HF_2^{\red}$ of dimension $2$, the pairing vanishing
on the $\RR$-summand.
The relative invariants will be
$\p^{w}(\X_1,e^{t\D_1}) \in V_2[[t]]$. Again we do not consider the odd
part, and we also ignore the extra $\RR$-summand.
\end{itemize}

\begin{prop}
\label{prop:concl}

\begin{enumerate}
  \item There is a vector space $V_4$ of dimension $4$ endowed with a
    symmetric bilinear form such that for every permissible $(X,\S)$
    and $\D \subset \X$ with $\bd \D =\SS^1$, we have
    $\p^w(\X,e^{t\D}) \in V_4[[t]]$. For $X=\X_1 \cup_Y \X_2$, 
    $D=\D_1 +\D_2$, $\bd \D_1 =-\bd \D_2 =\SS^1$, we have
    $$ D_X^{(w,\S)}(e^{tD})=
     <\p^{w_1}(\X_1,e^{t\D_1}),\p^{w_2}(\X_2,e^{t\D_2})>. $$
  \item There is a vector space $V_2$ of dimension $2$ endowed with a
    symmetric bilinear form such that for every permissible $(X,\S)$
    and $\D \subset \X$ with $\bd \D =\g$ not representing in homology
    a multiple of $[\SS^1]$, we have
    $\p^w(\X,e^{t\D}) \in V_2[[t]]$. For $X=\X_1 \cup_Y \X_2$, 
    $D=\D_1 +\D_2$, $\bd \D_1 =-\bd \D_2 =\g$, we have
    $$ D_X^{(w,\S)}(e^{tD})=
     <\p^{w_1}(\X_1,e^{t\D_1}),\p^{w_2}(\X_2,e^{t\D_2})>. $$
\end{enumerate}
\end{prop}

\section{Proof of Theorems}

\noindent {\em Proof of Theorem~\ref{thm:main1}.}

The fact that $X$ is of simple type will be proved in the proof
of theorem~\ref{thm:main3}.
Let us analyse the following list of examples 
(we use proposition~\ref{prop:DwS}
for finding the invariants $\Dws_X$).
  
\begin{itemize}
  \item $X$ a $K3$ surface blown-up twice with $E_1$ and $E_2$ the two
    exceptional divisors, $\S = S-E_1 -E_2$ for $S$ a
    tight surface of genus $2$ in $K3$, $w=E_1$, $D$ a cohomology class
    coming from the $K3$ such that $D \cdot S =1$, $D^2=0$. 
    We get $\Dws_X(e^{s\S+tD})=- e^{ts}\,{e^{2s}-e^{-2s} \over 4}$.
  \item $X$, $\S$, $D$ as before, but now $w \in H^2(K3)$, with $w
    \cdot S =1$. We will get $\Dws_X(e^{s\S+tD})= (-1)^{w^2 \over
    2}e^{ts}\,{e^{2s}+e^{-2s} \over 4} -{1\over 2}e^{-ts}$.
  \item $X$ a $K3$ surface, $\S$ a tight torus with an added trivial
    handle to make it of genus $2$, $w \in H^2(X ;\ZZ)$ such that $w
    \cdot \S=1$ and $D$ with $D\cdot \S=1$, $D^2=0$. Then 
    $\Dws_X(e^{s\S+tD})= -e^{-ts}$.
  \item $S = \CP^1 \x \S$, $w=\PD [\CP^1]$, $D=\CP^1$. Then
    $g(t,s)=\Dws_S(e^{s\S+tD})$ is a non-zero function with monomials of 
    degree at least three (since the smallest moduli space has dimension
    six).
\end{itemize}

We conclude that there are at least four functions, 
say $f_1=e^{ts+2s}$,$f_2=e^{ts-2s}$, $f_3=e^{-ts}$ and $f_4=g(t,s)$,
appearing in 
some $\Dws_X(e^{s\S+tD})$ (for different permissible pairs $(X,\S)$ and $D 
\cdot \S=1$), 
and linearly independent over $\cF(t)$, the 
field of (formal) Laurent series on $t$. We have the map
\begin{equation}
\label{eqn:alpha}
   < \cdot , \p^w(A,e^{t\De+s\S})>: V_4[[t]] \ar \RR^4[[t]]
\end{equation}
which assigns to $\p(t) \in V_4[[t]]$ a four-vector whose
i-th coordinate (actually we should tensor $V_4[[t]]$ and
$\RR^4[[t]]$ with $\cF(t)$, but we will not be explicit about this point)
is the coefficient (in $\cF(t)$) of
$f_i$ in $<\p(t),\p^w(A,e^{t\De+s\S})>$ (where $\De=\pt \x D^2 \subset A$).
Therefore $\p^w(\X,e^{t\D})$ is sent to $(c_{X,i}(t))$, the 
coefficients of $f_i$ in
$\Dws_X(e^{s\S+tD})$, where $D=\D +\De$ 
(so $\Dws_X(e^{s\S+tD}) = \sum c_{X,i}(t) \cdot 
f_i(t,s)$).

From the examples, the map above 
is an isomorphism (over $\cF(t)$), so we can push the product from 
$V_4[[t]]$ to $\RR^4[[t]]$ and we shall have a universal symmetric
matrix $M(t)=(M_{ij}(t))$ such that
$$
  \Dws_X(e^{tD})=\sum_{i,j} c_{X_1,i}(t)M_{ij}(t)c_{X_2,j}(t).
$$

\begin{rem}
\label{rem:extra}
Since the map~\eqref{eqn:alpha} is an isomorphism, $\Dws_X((x^2-4) 
e^{tD+s\S})=0$ if and only if
$\p^w(\X,(x^2-4)e^{t\D})=0$. 
\end{rem}

The image of all possible 
$\p^w(\X,e^{t\D})$ with $X$ of simple type, $b_1=0$ and $b^+>1$, is exactly
the three-dimensional subspace given by equating the last coordinate to
zero. So when $X_i$ are both of simple type with $b_1=0$ and $b^+>1$,
write $\DD^w_{X_1}(e^{\a})=e^{Q(\a)/2}\sum a_{i,w} e^{K_i\cdot \a}$ and
$\DD^w_{X_2}(e^{\a})=e^{Q(\a)/2}\sum b_{j,w} e^{L_j\cdot \a}$. 
Then 
  \begin{equation} 
    \left\{ \begin{array}{l}
     c_{X_1,1}(t) = e^{Q(tD_1)/2}
     \sum\limits_{K_j \cdot \S=2} \hspace{-3mm}
     a_{j,w}e^{tK_j \cdot D_1} \\
     c_{X_1,2}(t) = e^{Q(tD_1)/2}
     \sum\limits_{K_j \cdot \S=-2} \hspace{-3mm}
     a_{j,w}e^{tK_j \cdot D_1} \\
     c_{X_1,3}(t) = e^{-Q(tD_1)/2}
     \sum\limits_{K_j \cdot \S=0} \hspace{-3mm}         
     i^{-d_0} a_{j,w}e^{ti\,K_j \cdot D_1} \\
    c_{X_1,4}(t) = 0
    \end{array} \right. \label{eqn:coeff}
  \end{equation}
and 
  \begin{equation}
  \Dws_X(e^{tD})=\sum_{1 \leq i,j \leq 3}
  c_{X_1,i}(t)M_{ij}(t)c_{X_2,j}(t).
  \label{eqn:prod}\end{equation}

This expression is valid for any $D \in H_2(X)$ with $D|_Y=[\SS^1]$.
For $D \in \cH$ we have
$$  
\Dws_X(e^{tD})=\sum_{1 \leq i,j \leq 3}
  c_{X_1,i}(t)M_{ij}(t (D\cdot \S))c_{X_2,j}(t).
$$
Considering $D$, $D_1+r\S$, 
$D_2-r \S$ in \eqref{eqn:prod},
we get that $M_{ij}(t)=0$ for $i \neq j$, $1 \leq i,j \leq 3$.
Now consider the case in which both $X_i$
and $\S_i$ are as in the third example of the list. Then $X = X_1 \#_{\S}
X_2$ splits off a ${\Bbb S}^2 \x {\Bbb S}^2$, so its invariants are
zero. Therefore $M_{33}(t)=0$. So finally we have (using also
$D^2=D_1^2 +D_2^2$),
$$  
\Dws_X(e^{tD})= e^{Q(tD)/2}(\hspace{-5mm}\sum_{K_i \cdot \S=L_j \cdot
    \S = 2}\hspace{-5mm} M_{11}(t(D\cdot \S)) 
    a_{i,w}b_{j,w} \, e^{(K_i \cdot D_1 +  L_j
    \cdot D_2)t } +
$$
$$
  \hspace{-5mm}\sum_{K_i \cdot
    \S=L_j \cdot \S= -2}\hspace{-5mm} M_{11}(t(D\cdot \S)) 
    a_{i,w}b_{j,w} \, e^{(K_i \cdot
    D_1 +  L_j \cdot D_2)t}).
$$

  Let us now compute $M_{11}(t)$ and $M_{22}(t)$. 
  By the universality and since all the manifolds involved can be chosen 
  of simple type, one has $M_{11}(t) = \sum c_n\, e^{nt}$ 
  and $M_{22}(t)= \sum d_n\, e^{nt}$, finite sums of exponentials.
  Let $S=\CP^2 \# 10\overline{\CP}^2$ be the rational elliptic surface
  blown-up once. Denote by $E_1, \dots, E_{10}$ the exceptional
  divisors and let $T_1=C-E_1- \cdots -E_9$, $T_2=C-E_1- \cdots -E_8-
  E_{10}$, where $C$ is the cubic curve in $\CP^2$. So $T_1$ and $T_2$
  can be represented by smooth tori of self-intersection zero and with
  $T_1 \cdot T_2 =1$. We can glue two copies of $S$ along $T_1$. The
  result is a K3 surface $S \#_{T_1} S$ 
  blown-up twice. The $T_2$ pieces glue
  together to give a genus $2$ Riemann surface $\S_2$ of
  self-intersection
  zero which intersects $T_1$ in one point. This is actually the  
  pair $(B,\S)$ we introduced before the
  statement of theorem~\ref{thm:main3}.
  Now set $X= (S \#_{T_1} S ) \#_{\S_2} (S \#_{T_1} S )$, which
  is of simple type (by~\cite{KM}, since it contains a torus of  
  self-intersection $0$ 
  intersecting an embedded $(-2)$-sphere transversely in one point).
  Now call
  $\S=\S_2$ and get $D$ piecing together both $T_1$'s in $S \#_{T_1} S$.
  So (choose $w=T_1$ on $S \#_{T_1} S$) 
  $$
    D^{(D,\S)}_X(e^{tD+s\S}) = e^{Q(tD+s\S)/2}(\hspace{-4mm} \sum_{K_i
    \cdot \S=L_j \cdot 
    \S = 2} \hspace{-4mm} c_na_ib_j \, e^{2s+nt} +\hspace{-4mm}
    \sum_{K_i \cdot \S=L_j \cdot \S= -2}\hspace{-4mm} d_na_ib_j \,
    e^{-2s+nt}) = 
  $$
  $$ 
   =e^{ts} (\sum {c_n \over 16} e^{2s+nt} +\sum {d_n \over 16}
   e^{-2s+nt}),
  $$
  since $T_1$ evaluates $0$ on basic classes being a torus of
  self-intersection zero (the coefficient $1 \over 16$ appears from
  the explicit computation of the basic classes of the K3 surface
  blown-up in two points, see below~\eqref{eqn::}).
  The trick is now to use the symmetry fact that $X= (S \#_{T_2} S
  ) \#_{\S_1}  (S \#_{T_2} S )$, where $\S_1$ comes from gluing
  together both $T_1$'s. Under this diffeomorphism $D=\S_1$ and 
  $\S$ comes from piecing together both $T_2$'s in $S \#_{T_2} S$. Hence
  $$
     D^{(\S,D)}_X(e^{tD+s\S}) = e^{ts} (\sum {c_n \over 16} e^{2t+ns}
     +\sum {d_n \over 16}   e^{-2t+ns}).
  $$
  Both expressions are equal, and equal to $\DD_X^{D+\S}(e^{tD+s\S})$.
  From here we deduce that $c_n=0$ unless $n=\pm 2$ and $d_n=0$ unless
  $n= \pm 2$. Also $c_{-2}=d_2$. Put $l=c_2+c_{-2}$, so  
  $D^{(D,\S)}_X(e^{s\S}) = {l \over 16} e^{2s} - {l \over 16} e^{-2s}$
  (note that $d_0(X,D)=-15$ is odd).
  So $c_2 - d_{-2} =2\,l$. But $c_2 = \pm d_{-2}$, so it has to be
  $c_{-2}=d_2=0$ and $c_2=-d_{-2}=l$. Thus
  \begin{equation}\label{eqn:***}
     D^{(D,\S)}_X(e^{tD+s\S}) = e^{ts} ({l \over 16} e^{2s+2t}
     - {l \over 16}   e^{-2s-2t}).
  \end{equation}
  So $M_{11}(t)=l \, e^{2t}$ and $M_{22}(t)=-l \, e^{-2t}$. 
  To get the theorem it only remains to prove
  
\begin{lem}
\label{lem:l=32}
  $l=-32.$
\end{lem}

\begin{pf}
  Let $e_i =\p^w(A,\S^i) \in HF_*(Y)$, $i=0,1,2,3$. Then $\{e_i\}$ is a
  basis for $HF_{\mathrm{even}}(Y)$, since the latter is a vector
  space of dimension $4$ and the intersection matrix for $(e_i\cdot e_j)$
  is invertible. Actually, it is
  $$ 
    N =\left( \begin{array}{cccc} 0 & 0& 0& -1/2 \\  0& 0& -1/2  & 0\\
     0& -1/2 &0 & -2 \\ -1/2 & 0& -2 &0 \end{array} \right). 
  $$
  To check this we note that $e_i \cdot e_j=<\p^w(A,\S^i),\p^w(A,\S^j)>=
  \Dws_S(\S^{i+j})$, $S=\S\x\CP^1=A \cup_Y A$, $w=\CP^1$, so we 
  only need to find $\Dws_S(\S^3)$ and $\Dws_S(\S^5)$. 
  For the first one, the moduli space is $\M$, which is
  six-dimensional. Then $\mu(\S)^3 = 1/2$, with $\mu(\S) \in H^2(\M)$, 
  from~\cite{Thaddeus}. This invariant is computed using the complex
  orientation of the moduli space which differs from the one we use 
  by a factor $\f_S(w)=(-1)^{K_S \cdot w +w^2 \over 2}=-1$. So
  $\Dws_S(\S^3) =-1/2$. For the 
  second one, the moduli space is ten-dimensional, corresponding to
  $w=\CP^1 + \S$ and polarisation close to $\S$. 
  For a polarisation close to $\CP^1$, the moduli space is 
  empty~\cite{Qin}. There 
  is only one wall corresponding to $\z=-\S+\CP^1$. 
  Now we can apply the
  formulas in~\cite{papervic} for wall-crossing when the irregularity is
  not zero, noting that $\z$ is a good wall. This gives $\Dws_S(\S^5)=
  -2$. Also this can be computed directly with an explicit description
  of the algebraic moduli space and we propose this calculation as
  a good exercise.
  
  Now consider the pair $(B,\S)$. Recall that $B$ is the $K3$ surface 
  blown-up in two points. Let $E_1$ and $E_2$ be the two exceptional
  divisors. We have $\DD_B(e^{\a})=
  e^{Q(\a)/2} \sinh (E_1 \cdot  \a) \sinh (E_2 \cdot \a)$, so by 
  proposition~\ref{prop:DwS} with $w=T_1$,
\begin{equation}
     \Dws_B(e^{\a})= e^{Q(\a)/2}{1 \over 2} \cosh ((E_1+E_2) \cdot 
     \a) +e^{-Q(\a)/2} {1 \over 2} \cos ((E_1-E_2) \cdot  \a).
\label{eqn::}\end{equation}
  Then $\Dws_B(e^{s\S})={1 \over 2} \cosh(2s)$, so in the  basis 
  dual to $\{e_i\}$, $\p^w(\B,1) =(1/2,0,2,0)$   
  and $\p^w(\B,\S) =(0,2,0,8)$. Therefore for $C=B \#_{\S} B$,
  $$
    \Dws_C(\S)=<\p^w(\B,1), \p^w(\B,\S)>= (1/2,0,2,0) N^{-1}\,N\,N^{-1} 
    {\left( \begin{array}{c}
    0 \\ 2\\0\\8 \end{array}\right)}=-8.
  $$
  From formula~\eqref{eqn:***}, we get $\Dws_C(\S)={l \over 4}$, so
  $l=-32$.
\end{pf}

\noindent {\em Proof of Theorem~\ref{thm:main2}.}

  We are going to check that $\Dws_X((x^2-4)^2 e^{tD})=0$, for all $D \in
  H_2(X)$ with $D \cdot \S =1$. This is clearly enough to infer the
  result. 
  Put $D=\D +\De$. If $X$ is of simple type with $b_1=0$ and $b^+>1$, 
  we have
  $$
    0=\Dws_X((x^2-4)e^{tD+s\S})=
     <\p^w(\X,e^{t\D}),\p^w(A,(x^2-4)e^{t\De+s\S})>.
  $$
  The vectors $\p^w(\X,e^{t\D})$ (with $X$ being of simple
  type with $b_1=0$ and $b^+>1$) generate a the $3$-dimensional 
  subspace $V_3[[t]]$
  in $V_4[[t]]$ given by equating the last coordinate to zero.
  Then $\p^w(A,(x^2-4)e^{t\De+s\S})$ lies in the subspace orthogonal
  to $V_3[[t]]$. As the pairing in $V_4[[t]]$ is non-degenerate and 
  $V_3[[t]]$ contains an isotropic vector (from the computation of 
  the $M_{ij}(t)$ in the proof of theorem~\ref{thm:main1}, the
  intersection matrix restricted to $V_3[[t]]$ is degenerate),
  $\p^w(A,(x^2-4)e^{t\De+s\S})$ is isotropic and hence
  $$
     <\p^w(A,(x^2-4)e^{t\De+s\S}), \p^w(A,(x^2-4)e^{t\De})>= \Dws_{\S \x
     \CP^1} ((x^2-4)^2e^{t\CP^1+s\S})=0,
  $$
  from where $\p^w(A,(x^2-4)^2e^{t\De})=0$ (remark~\ref{rem:extra})
  and hence the result.

\noindent {\em Proof of Theorem~\ref{thm:main3}.}

Recall the permissible pair $(B,\S)$, where $B$ is
the $K3$ surface blown-up in two points with $E_1$ and $E_2$
the exceptional divisors, and $\S=S-E_1-E_2$ is the proper transform of 
a tight embedded surface $S\subset K3$ of genus $2$.
Call $C=B \#_{\S} B$ the double of $B$, i.e.
the connected sum of $B$ with itself with the identification
which is given by the natural orientation reversing diffeomorphism of
$Y=\bd B^o$ to itself. 
As in the proof of theorem~\ref{thm:main1},
we choose $D \subset C$ to be the embedded surface obtained 
by piecing together
two fibres of the natural elliptic fibration of $B$. Then $D$ is
a genus $2$ Riemann surface of self-intersection zero. Also take
$w=\PD [D] \in H^2(X; \ZZ)$. Then equation~\eqref{eqn:***} gives
$$
    \Dws_C(e^{tD+s\S})=-e^{ts}(2\, e^{2s+2t} -2\,e^{-2s-2t}).
$$ 
We can take a collection $\a_i$, $1 \leq i \leq 4$, of loops 
in a fibre $\S \subset \bd B^o$, which together with $\SS^1$
form a basis for $H_1(Y)$, such that they can be
capped off with embedded $(-1)$-discs $D_i$ (writing $B = S \#_{T_1}
S$, as in the proof of theorem~\ref{thm:main1}, we consider
the vanishing discs of the elliptic fibration of $S$ with fibre 
$T_2$, see~\cite[page 167]{FM}, since they do not intersect $T_1$).
Now these discs can
be glued together pairwise when forming $C=\B\cup_Y \B$ to give a
collection of $(-2)$-embedded spheres $S_i=D_i \cup_{\a_i} D_i$.
Every one of these discs has a dual torus $T_i$, by considering another
loop in $\S \subset \bd \B$, say $\b_i$, with $\a_i \cdot
\b_i =1$, and putting $T_i =\b_i \x \SS^1 \subset \Y$. Then the
elements $S_i +T_i$ are represented by embedded tori of
self-intersection zero. Hence the manifold $C$ is of simple
type~\cite{KM}, and the
basic classes evaluate zero on $T_i$ and on $S_i +T_i$. Our conclusion
is
$$ 
    \Dws_C(e^{\a})=-4\,e^{Q(\a)/2}\sinh (K \cdot \a),
$$
with $K \in H^2(C;\ZZ)$ being the only cohomology class with
  \begin{itemize}
  \item $K \cdot \a =(E_1 +E_2) \cdot \a$ for $\a \in H_2(B^o)$.
  \item $K \cdot \S=K \cdot D =2$.
  \item $K \cdot S_i= K \cdot T_i =0$, for all $i$.
  \end{itemize}
We split $K$ into two symmetric pieces $K^o \subset \B$. The boundary of
$K^o$ is $\bd K^o =2\SS^1$ and $(K^o)^2 =2$ since $K^2=4$.

Analogously, the manifold $C_2=C \#_{\S} B$ is of simple type and 
$$ 
    \Dws_{C_2}(e^{\a})=32\,e^{Q(\a)/2}\cosh (K_2 \cdot \a).
$$
for a unique $K_2 \in H^2(C_2; \ZZ)$.
Let $D_2$ be obtained gluing the $D$ coming from $C$ with 
one fibre of the elliptic fibration of $B$. Then
$$
    \Dws_{C_2}(e^{tD_2+s\S})=e^{ts}(16\, e^{2s+2t} + 16\,e^{-2s-2t}).
$$ 

So there are two functions, $\tilde f_1=e^{ts+2s}$,$\tilde f_2=e^{ts-2s}$,
appearing in 
some $\Dws_{\tilde X}(e^{s\S+tD})$ (for different permissible
$(X,\S)$ with $b_1=0$, $\tilde X= \X \cup_Y B^o$,
$D \in H_2(\tilde X)$, $D=\D+\D_B$, $\bd \D=\g$), 
and linearly independent over $\cF(t)$.
Now we mimic the reasoning of the proof of theorem~\ref{thm:main1}. We 
have a map
$$
   < \cdot , \p^w(\B,e^{t\D_B+s\S})>: V_2[[t]] \ar \RR^2[[t]]
$$
which is an isomorphism (over $\cF(t)$), such that
$\Dws_{\tilde X}(e^{s\S+tD})= \sum c_{\tilde X,i}(t) \cdot 
\tilde f_i(t,s)$. We shall have a universal symmetric
matrix $\tilde M(t)=(\tilde M_{ij}(t))$ such that
  $$
  \Dws_{X}(e^{tD})=\sum_{i,j} c_{\tilde X_1,i}(t)\tilde M_{ij}
  (t)c_{\tilde X_2,j}(t).
  $$ 
This expression is valid for any $D \in H_2(X)$ with $D|_Y =[\g]\in H_1(Y)$, 
$\g \subset \S \subset Y$ an embedded curve.
Now $\p^w(\B,(x^2-4)e^{t\D_B})=0$, since $\Dws_C((x^2-4)e^{tD+s\S})=0$,
as $C$ is of simple type.
Therefore $\tilde X_i=\X_i \cup_Y \B$ are of simple type. Also this
implies that
$\p^w(\X_i,(x^2-4)e^{t\D})=0$ and hence that $X=\X_1 \cup_Y \X_2$ is of
simple type. Now
  \begin{equation}
  \left\{ \begin{array}{l}
     c_{\tilde X_1,1}(t) = e^{Q(tD_1)/2}
     \sum\limits_{\tilde K_i \cdot \S=2} \hspace{-3mm}
     \tilde a_{i,w}e^{t\tilde K_i \cdot D_1} \\
     c_{\tilde X_1,2}(t) = e^{Q(tD_1)/2}
     \sum\limits_{\tilde K_i \cdot \S=-2} \hspace{-3mm}
     \tilde a_{i,w}e^{t\tilde K_i \cdot D_1} \\
    \end{array} \right. \label{eqn:coeff2}
  \end{equation}
  
Again, as in the proof of theorem~\ref{thm:main1},
we get that $\tilde M_{ij}(t)=0$ for $i \neq j$.
Consider $X_1=B$, $X_2=B$, $X=C$, $\tilde X_1=C$, $\tilde X_2=C$, 
let $D_C=\D_B+\D_B \in H_2(C)$ and put $D_1=D_2=D_C$, $D=D_C$.
We get that $\tilde M_{11}(t) = -{1 \over 2}
e^{Q(tD_C)/2}$, 
$\tilde M_{22}(t) = {1 \over 2}
e^{Q(tD_C)/2}$. Now we have that for any $X=\X_1 \cup_Y \X_2$ and
$D=\D_1 +\D_2$ with $\bd \D_1=-\bd \D_2=\g$, it is
$Q(tD_C)/2= Q(tD)/2 -Q(tD_1)/2-Q(tD_2)/2$, where $D_i=\D_i +\D_B \in
H_2(\tilde X_i)$.
From this we get the sought expression in the statement of 
theorem~\ref{thm:main3}, for any $D \in H_2(X)$ with $D|_Y \in H_1(Y)$ 
satisfying that $p_*(D|_Y) \in H_1(\S)$ ($p: Y \to \S$ the projection)
is primitive and non-zero.
Since we have chosen the map $D \mapsto (D_1, D_2)$ to be linear, 
this finishes the proof.

\end{document}